\documentclass{amsart}
\usepackage{amsmath}
\usepackage{amstext}
\usepackage{amssymb}
\usepackage{amsthm}
\swapnumbers
\theoremstyle{plain}
\newtheorem{thm}{Theorem}[section]
\newtheorem{lem}[thm]{Lemma}
\newtheorem{prop}[thm]{Proposition}
\newtheorem{cor}[thm]{Corollary}

\theoremstyle{definition}

\newtheorem{ntn}[thm]{Notation}

\theoremstyle{remark}
\newtheorem*{note}{Note}
\newtheorem{rmk}[thm]{Remark}

 \DeclareMathOperator{\height}{ht}
 \DeclareMathOperator{\lm}{lm}
\DeclareMathOperator{\lc}{lc} \DeclareMathOperator{\lcm}{lcm}
\DeclareMathOperator{\lt}{lt} \DeclareMathOperator{\nd}{end}
 
\DeclareMathOperator{\beg}{beg} 
\DeclareMathOperator{\pd}{proj\,dim}
\DeclareMathOperator{\Ker}{Ker}
\DeclareMathOperator{\Coker}{Coker} \DeclareMathOperator{\Ima}{Im}
\DeclareMathOperator{\depth}{depth}

\DeclareMathOperator{\Supp}{Supp} 
 
\DeclareMathOperator{\Ass}{Ass} 
 
\DeclareMathOperator{\Var}{Var} 
\DeclareMathOperator{\Spec}{Spec}
\DeclareMathOperator{\Proj}{Proj}

\def\Z{\mathbb Z}
\def\N{\mathbb N}

\def\fa{{\mathfrak{a}}}

\def\fm{{\mathfrak{m}}}

\def\p{{\mathfrak{p}}}
\def\fp{{\mathfrak{p}}}

\def\fq{{\mathfrak{q}}}

\def\nn{\relax\ifmmode{\mathbb N_{0}}\else$\mathbb N_{0}$\fi}
\def\lra{\longrightarrow}

\begin{document}

\title[ASSOCIATED PRIMES OF GRADED COMPONENTS OF LOCAL COHOMOLOGY
MODULES]{ASSOCIATED PRIMES OF GRADED COMPONENTS OF LOCAL
COHOMOLOGY MODULES}
\author{MARKUS P. BRODMANN}
\address[Brodmann]{Institut f\"{u}r Mathematik, Universit\"{a}t Z\"{u}rich,
Winterthurerstrasse
190, 8057 Z\"{u}rich, Switzerland\\
{\it Fax number}: 0041-1-635-5706} \email{Brodmann@math.unizh.ch}
\author{MORDECHAI KATZMAN}
\address[Katzman]{Department of Pure Mathematics,
University of Sheffield, Hicks Building, Sheffield S3 7RH, United Kingdom\\
{\it Fax number}: 0044-114-222-3769}
\email{M.Katzman@sheffield.ac.uk}
\author{RODNEY Y. SHARP}
\address[Sharp]{Department of Pure Mathematics,
University of Sheffield, Hicks Building, Sheffield S3 7RH, United Kingdom\\
{\it Fax number}: 0044-114-222-3769}
\email{R.Y.Sharp@sheffield.ac.uk}
\thanks{The third author was partially supported by the Swiss National
Foundation (Project number 20-52762.97).}


\date{\today}


\begin{abstract}
The $i$-th local cohomology module of a finitely generated graded
module $M$ over a standard positively graded commutative
Noetherian ring $R$, with respect to the irrelevant ideal $R_+$,
is itself graded; all its graded components are finitely generated
modules over $R_0$, the component of $R$ of degree $0$.  It is
known that the $n$-th component $H^i_{R_+}(M)_n$ of this local
cohomology module $H^i_{R_+}(M)$ is zero for all $n>> 0$. This
paper is concerned with the asymptotic behaviour of
$\Ass_{R_0}(H^i_{R_+}(M)_n)$ as $n \rightarrow -\infty$.

The smallest $i$ for which such study is interesting is the
finiteness dimension $f$ of $M$ relative to $R_+$, defined as the
least integer $j$ for which $H^j_{R_+}(M)$ is not finitely
generated.  Brodmann and Hellus have shown that
$\Ass_{R_0}(H^f_{R_+}(M)_n)$ is constant for all $n < < 0$ (that
is, in their terminology, $\Ass_{R_0}(H^f_{R_+}(M)_n)$ is
asymptotically stable for $n \rightarrow -\infty$). The first main aim
of this paper is to identify the ultimate constant value (under the mild
assumption that $R$ is a homomorphic image of a regular ring): our answer
is precisely the set of contractions to $R_0$ of certain relevant primes
of $R$ whose existence is confirmed by Grothendieck's Finiteness Theorem
for local cohomology.

Brodmann and Hellus raised various questions about such asymptotic
behaviour when $i > f$. They noted that Singh's study of a
particular example (in which $f = 2$) shows that
$\Ass_{R_0}(H^3_{R_+}(R)_n)$ need not be asymptotically stable for
$n \rightarrow -\infty$. The second main aim of this paper is to
determine, for Singh's example, $\Ass_{R_0}(H^3_{R_+}(R)_n)$ quite
precisely for every integer $n$, and, thereby, answer one of the
questions raised by Brodmann and Hellus.
\end{abstract}

\maketitle

\setcounter{section}{-1}
\section{\bf Introduction}
\label{in}

Let $R = \bigoplus_{n \in \nn} R_n$ be a positively graded
commutative Noetherian ring which is standard in the sense that $R
= R_0[R_1]$, and set $R_+ := \bigoplus_{n \in \N} R_n$, the
irrelevant ideal of $R$. (Here, $\nn$ and $\N$ denote the set of
non-negative and positive integers respectively; $\Z$ will denote
the set of all integers.) Let $M = \bigoplus_{n \in \Z} M_n$ be a
non-zero finitely generated graded $R$-module. This paper is
concerned with the behaviour of the graded components of the
graded local cohomology modules $H^i_{R_+}(M)~(i \in \nn)$ of $M$
with respect to $R_+$.

It is known (see \cite[15.1.5]{BS}) that there exists $r \in \Z$ such that
$H^i_{R_+}(M)_n = 0$ for all $i \in \nn$ and all $n \geq r$, and that
$H^i_{R_+}(M)_n$ is a finitely generated $R_0$-module for all $i \in \nn$
and all $n \in \Z$.
Set
$$
f := f_{R_+}(M) =
\inf \left\{ i \in \N : H^i_{R_+}(M)
\mbox{~is not finitely generated} \right\},
$$
the finiteness dimension of $M$ relative to $R_+$: see {\rm
\cite[9.1.3]{BS}}. We assume that $f$ is finite.  M. Brodmann and
M. Hellus have shown in \cite[Proposition 5.6]{BroHel00} that
$\Ass_{R_0}(H^f_{R_+}(M)_{n})$ is constant for all $n < < 0$. The
first part (\S \ref{ab}) of this paper
determines the ultimate constant value under the mild restriction that
$R$ is a homomorphic image of
a regular (commutative Noetherian) ring; the main
result is related to Grothendieck's Finiteness Theorem for local
cohomology, which (under the specified restriction) gives an alternative
description of $f$.   Let $\mbox{\rm *}\Spec (R)$ denote the set
of graded prime ideals of $R$, and $\Proj (R)$ denote the set $\{
\p \in \mbox{\rm *}\Spec (R) : \p \not\supseteq R_+ \}$. Write
$$
\lambda^{R_+}_{R_+}(M) := \inf \left\{ \depth_{R_{\mathfrak{p}}}
M_{\mathfrak{p}} + \height (R_+ + {\mathfrak{p}})/{\mathfrak{p}} :
{\mathfrak{p}} \in \Proj (R) \right\}.
$$
(We interpret the depth of a zero module as $\infty$.) It is a
consequence of Grothendieck's Finiteness Theorem
\cite[Expos\'e VIII, Corollaire 2.3]{Groth68} that, when $R$ is
a homomorphic image of a regular ring,
$$
f = \lambda^{R_+}_{R_+}(M) = \inf \left\{
\depth_{R_{\mathfrak{p}}} M_{\mathfrak{p}} + \height (R_+ +
{\mathfrak{p}})/{\mathfrak{p}} : {\mathfrak{p}} \in \Proj (R)
\right\}.
$$
(See {\rm \cite[13.1.17]{BS}.}) The main result of \S \ref{ab} is
that, under the assumption that $R$ is a homomorphic image of a
regular ring,
$$
\left\{ \fp \cap R_0 : \fp \in \Proj (R) \text{~and~}
\depth_{R_{\mathfrak{p}}}
M_{\fp} + \height (\fp + R_+)/\fp = f\right\} =
\Ass_{R_0}(H^f_{R_+}(M)_n) \quad \text{~for all~} n < < 0.
$$

The final \S \ref{se} is concerned with the asymptotic behaviour
of $\Ass_{R_0}(H^i_{R_+}(M)_{n})$ as $n \rightarrow -\infty$ when
$i
> f$. Brodmann and Hellus say that $\Ass_{R_0}(H^i_{R_+}(M)_n)$ is
{\em asymptotically stable\/} (respectively {\em asymptotically
increasing\/}) for $n \rightarrow -\infty$ if there exists $n_0
\in \Z$ such that $\Ass_{R_0}(H^i_{R_+}(M)_n) =
\Ass_{R_0}(H^i_{R_+}(M)_{n_0})$ (respectively
$\Ass_{R_0}(H^i_{R_+}(M)_n )\subseteq
\Ass_{R_0}(H^i_{R_+}(M)_{n-1})$) for all $n \leq n_0$. They used
an example of A. Singh \cite[\S 4]{Singh00} to show that, when $i
> f$, $\Ass_{R_0}(H^i_{R_+}(M)_n)$ need not be asymptotically
stable for $n \rightarrow -\infty$. In \S \ref{se}, we use
Gr\"obner basis techniques to show that, for Singh's example
$$
R = \Z[X,Y,Z,U,V,W]/(XU + YV + ZW),
$$
where the polynomial ring $\Z[X,Y,Z,U,V,W]$ is graded so that its
$0$-th component is $\Z[X,Y,Z]$ and $U,V,W$ have degree $1$, we
have
$$
\Ass_{R_0}(H^3_{R_+}(R)_{-d}) = \left\{(X,Y,Z)\right\} \cup
\left\{(p,X,Y,Z) : p \in \Pi (d-2) \right\} \quad \mbox{~for all~}
d \geq 3,
$$
where $$\Pi (d-2) := \left\{ p : p \mbox{~is a prime factor of~}
{{d-2} \choose {i}} \mbox{~for some~} i \in \{0, \dots,
d-2\}\right\}.$$ It follows that $\Ass_{R_0}(H^3_{R_+}(R)_{n})$ is
not asymptotically increasing for $n \rightarrow -\infty$, and
this settles a question raised by Brodmann and Hellus.

\section{\bf Asymptotic behaviour at the finiteness dimension}
\label{ab}

\begin{ntn}\label{mi.0} {\rm The notation introduced in the above
\S \ref{in} will be maintained for the whole paper. We shall only
assume that $R$ is a homomorphic image of a regular ring when this
is explicitly stated. Here we
introduce additional notation.

We use,
for $j \in \Z$, the notation $L_j$ to denote the $j$-th component
of a $\Z$-graded module $L$, and $ (\: {\scriptscriptstyle
\bullet} \:)(j)$ to denote the $j$-th shift
functor on the category
of graded $R$-modules and homogeneous homomorphisms (by
`homogeneous' here, we mean `homogeneous of degree zero').
It will be convenient to have available
the concepts of the {\em end\/} and {\em beginning\/}
($\beg (L)$) of the graded $R$-module $L = \bigoplus
_{n \in \Z}L_n$, which are defined by $$\nd (L)
:= \sup \left\{ n \in \Z : L_n \neq 0\right\}
\quad \mbox{~and~} \quad \beg (L)
:= \inf \left\{ n \in \Z : L_n \neq 0\right\}.
$$
(Note that $\nd (L)$ could be $\infty$, and that the supremum of the empty set
of integers is to be taken as $- \infty$; similar comments
apply to $\beg (L)$.)

For $\fp \in \Spec (R)$, we abbreviate $\depth_{R_{\mathfrak{p}}}
M_{\mathfrak{p}}$ by $\depth M_{\mathfrak{p}}$ and the projective
dimension $\pd_{R_{\mathfrak{p}}} M_{\mathfrak{p}}$ by $\pd
M_{\mathfrak{p}}$.}
\end{ntn}

\begin{lem}\label{bo.14} (The notation is as in {\rm \S \ref{in}} and {\rm \ref{mi.0}}.)
Let $\fp \in \Proj (R) \cap \Ass _RM$ be such that
$\height (\fp + R_+)/\fp  = 1$.
Set $\fp_0 = \fp \cap R_0$.
Then $\fp_0 \in \Ass_{R_0}(H^1_{R_+}(M)_n)$ for all $n < \beg (M)$.
\end{lem}

\begin{proof} Set $\overline{M} := M/\Gamma_{R_+}(M)$, and note
that, by \cite[2.1.12 and 2.1.7(iii)]{BS}, $$\Ass_R(\overline{M}) =
\Proj(R) \cap \Ass_RM$$ and there is a homogeneous isomorphism
$H^1_{R_+}(M) \cong H^1_{R_+}(\overline{M})$. We therefore can, and do,
assume that $\Gamma_{R_+}(M) = 0$ in the remainder of this proof.

We now use homogeneous localization at $\fp + R_+$ to see that is
is enough to prove the claim under the additional hypotheses
that $R$ is $\mbox{\rm *}$local with unique $\mbox{\rm *}$maximal ideal
$\fm$, and that $\fm_0 := \fm \cap R_0 = \fp_0$.
The assumptions that $R$ is standard and $\mbox{\rm *}$local with
$\fm_0 = \fp_0$, and that
$\height (\fp + R_+)/\fp  = 1$, ensure that there exists
$g_1 \in R_1 \setminus \fp$, and that, then, $\sqrt{\fp + g_1R}
= \fp + R_+$.

Now there exists $t \in \Z$ such that $M$ has a graded $R$-submodule $N$
homogeneously isomorphic to $(R/\fp)(-t)$. We now consider the
ideal transform $D_{Rg_1}(N)$ of $N$ with respect to
$Rg_1$: this is naturally graded, and since $g_1$ is a non-zerodivisor on
$R/\fp$, the description of this ideal transform afforded by
\cite[Theorem 2.2.16]{BS} shows that (a) multiplication by $g_1$
provides a homogeneous isomorphism $D_{Rg_1}(N) \stackrel{\cong}{\lra}
D_{Rg_1}(N)(1)$, and that (b) $\fp_0 \in \Ass_{R_0}
((D_{Rg_1}(N)_n)$ for all $n \in \Z$.

Point (a) leads to the conclusion that multiplication by $g_1$
provides a homogeneous isomorphism $H^i_{R_+}(D_{Rg_1}(N))
\stackrel{\cong}{\lra}
H^i_{R_+}(D_{Rg_1}(N))(1)$ for each $i \in \nn$, and so, in particular,
for $i = 0$ and $1$. Since $g_1 \in R_+$, we conclude that
$H^i_{R_+}(D_{Rg_1}(N)) = 0$ for $i = 0,1$.

By \cite[2.2.4]{BS}, the natural (homogeneous) monomorphism
$\eta_N : N \lra D_{Rg_1}(N)$ has cokernel isomorphic to $H^1_{Rg_1}(N)$.
But, since $\sqrt{\fp + g_1R} = \fp + R_+$ and $N$ is
annihilated by $\fp$, there is a (homogeneous)
isomorphism $H^1_{Rg_1}(N)\cong H^1_{R_+}(N)$. Thus the cokernel of the
monomorphism $\eta_N$ is $R_+$-torsion.

It now follows from \cite[2.2.13 and 12.4.2(ii)]{BS} that there
is a homogeneous isomorphism $D_{Rg_1}(N) \cong D_{R_+}(N)$. We can thus
conclude that $\fp_0 \in \Ass_{R_0}
((D_{R_+}(N)_n)$ for all $n \in \Z$.
We now note that, since $D_{R_+}$ is a left exact functor, there is a
homogeneous $R$-monomorphism $D_{R_+}(N) \lra D_{R_+}(M)$; the result now
follows from the exact sequence
$$
0 \lra M \lra D_{R_+}(M) \lra H^1_{R_+}(M) \lra 0
$$
of graded $R$-modules and homogeneous $R$-homomorphisms.
\end{proof}

For part of the proof of our main result of this section, we shall
be able to reduce to the case where $R_0$ is a regular local ring
and $R = R_0[X_1, \ldots, X_r]$ is a polynomial ring over $R_0$ in which
the independent indeterminates $X_1, \ldots, X_r$ all have degree $1$.
This explains why several subsequent lemmas are concerned with this case.

\begin{lem}\label{bo.41} The notation is as in {\rm \S \ref{in}} and
{\rm \ref{mi.0}}. In addition, suppose that $(R_0,\fm_0)$ is a
regular local ring of dimension $d$ and that $R = R_0[X_1, \ldots,
X_r]$, a polynomial ring graded in the usual way. Suppose that
$\fp \in \Supp (M) \cap \Proj (R)$ is such that $\fp \cap R_0 =
\fm_0$. Then
$$
\depth M_{\fp} + \height (\fp + R_+)/\fp = d + r - \pd M_{\fp}.
$$
\end{lem}

\begin{proof} As $R$ is a catenary domain,
$$
\height (\fp + R_+)/\fp = \height (\fp + R_+) - \height \fp = d + r -
\height \fp.
$$
Moreover, by the Auslander-Buchsbaum-Serre Theorem,
$$
\depth M_{\fp} = \dim R_{\fp} - \pd M_{\fp} = \height \fp  - \pd M_{\fp}.
$$
\end{proof}

\begin{lem}\label{bo.42} The notation is as in {\rm \S \ref{in}} and {\rm \ref{mi.0}}.
In addition, suppose that $(R_0,\fm_0)$ is a regular local
ring of dimension $d$
and that $R = R_0[X_1, \ldots, X_r]$,
a polynomial ring graded in the usual way.

Let $(R'_0,\fm'_0)$ be a regular local
flat extension ring of $R_0$ such that $\fm_0R'_0 = \fm'_0$. Let
$R' = R \otimes_{R_0}R'_0$, which we identify with $R'_0[X_1, \ldots, X_r]$
in the obvious way. Let $M'$ denote the finitely generated graded
$R'$-module $M\otimes_RR'$, and let $\fp' \in \Proj (R')$ be such that
$\fp' \cap R'_0 = \fm'_0$. Set $\fp := \fp' \cap R$. Then $\fp \in \Proj (R)$
and $\fp \cap R_0 = \fm_0$; also
$$
\depth M_{\fp} + \height (\fp + R_+)/\fp \leq
\depth_{R'_{\fp'}} M'_{\fp'} + \height (\fp' + R'_+)/\fp'.
$$
\end{lem}

\begin{proof} Observe that there are $R'_{\fp'}$-isomorphisms
$$
M'_{\fp'} \cong (M\otimes_RR')\otimes_{R'} R'_{\fp'} \cong
M\otimes_R R'_{\fp'} \cong M_{\fp} \otimes_{R_{\fp}} R'_{\fp'}.
$$
As $R'_{\fp'}$ is a flat $R_{\fp}$-algebra, $\pd M'_{\fp'} \leq \pd M_{\fp}$.
Hence, by two uses of Lemma \ref{bo.41},
\begin{align*}
\depth M_{\fp} + \height (\fp + R_+)/\fp & =
d + r - \pd M_{\fp} \\
& \leq d + r - \pd M'_{\fp'} \\
& = \depth M'_{\fp'} + \height (\fp' + R'_+)/\fp'.
\end{align*}
\end{proof}

\begin{lem}\label{bo.43} The notation is as in {\rm \S \ref{in}} and
{\rm \ref{mi.0}}. In addition, suppose that $(R_0,\fm_0)$ is a
regular local ring of dimension $d$ such that the field
$R_0/\fm_0$ is algebraically closed and that $R = R_0[X_1, \ldots,
X_r]$, a polynomial ring graded in the usual way.

Suppose that $r > 1$, that $f_{R_+}(M) = r$ and that $\fm_0 \in
\Ass_{R_0}(H^r_{R_+}(M)_n)$ for all $n < < 0$. Then there exists $y \in
R_1 \setminus \fm_0R_1$ such that $y$
is a non-zerodivisor on $M/\Gamma_{R_+}(M)$,
that $f_{R_+}(M/yM) = r-1$ and that
$$
\fm_0 \in
\Ass_{R_0}(H^{r-1}_{R_+}(M/yM)_n) \quad \text{~for all~} n < < 0.
$$
\end{lem}

\begin{proof} Set $\overline{M} := M/\Gamma_{R_+}(M)$.
For a homogeneous element $y$ of $R$, we have
homogeneous isomorphisms
$$
\overline{M}/y\overline{M} \cong M/(yM + \Gamma_{R_+}(M)) \cong
(M/yM)/((yM + \Gamma_{R_+}(M))/yM),
$$
so that there are homogeneous isomorphisms $H^i_{R_+}(M) \cong
H^i_{R_+}(\overline{M})$ and $H^{i}_{R_+}(M/yM) \cong
H^{i}_{R_+}(\overline{M}/y\overline{M})$ for all
$i > 0$. We may therefore replace
$M$ by $\overline{M}$. We therefore assume that $\Gamma_{R_+}(M) = 0$
and $\Ass_RM \subseteq \Proj (R)$.

Now let $\fp \in \Ass_RM$ and set $\fp_0 := \fp \cap R_0$. Then, since
$R$ is a regular, and therefore catenary, domain,
\begin{align*}
\height_{R_0}\fp_0 + r - \height \fp & =
\height (\fp_0R + R_+) - \height \fp \\
& = \height (\fp_0R + R_+)/\fp = \height (\fp + R_+)/\fp \\
& = \depth M_{\fp} + \height (\fp + R_+)/\fp \\
& \geq f_{R_+}(M) = r.
\end{align*}
(We have used Grothendieck's Finiteness Theorem to obtain the inequality.)
Therefore $\height \fp_0R = \height_{R_0}\fp_0 = \height \fp$, so that
$\fp = \fp_0R$ and $\fp \subseteq \fm_0R$. It therefore follows that, if we
let $U$ denote the subset of $R_1 \setminus \fm_0R_1$ defined by
$$
U := \left\{ a_1X_1 + a_2X_2 : (a_1,a_2) \in R_0 \times R_0 \setminus
(\fm_0 \times \fm_0) \right\},
$$
then $U \cap \fp = \emptyset$. Therefore each element of $U$ is a
non-zerodivisor on $M$.

Set $J := \Gamma_{\fm_0R}(H^r_{R_+}(M)) = \bigoplus_{n \in \Z}
\Gamma_{\fm_0}(H^r_{R_+}(M)_n)$. The hypotheses ensure that $J$ is not a
finitely generated $R$-module. We shall show that one of the elements of
$U$ can be taken for $y$. To achieve this, we suppose that, for all $x \in U$,
there exists $n_x \in \Z$ such that, for all $n \leq n_x$, it is the
case that $\fm_0 \not\in \Ass_{R_0}(H^{r-1}_{R_+}(M/xM)_n)$, and we seek
a contradiction.

This supposition means that, for each $x \in U$, we have
$\Gamma_{\fm_0}(H^{r-1}_{R_+}(M/xM)_n) = 0$ for all $n \leq n_x$. Since
$f_{R_+}(M) = r$, there exists $\widetilde{n} \in \Z$
such that $H^{r-1}_{R_+}(M)_n
= 0$ for all $n \leq \widetilde{n}$. For each $x \in U$, the
application of local cohomology with respect to $R_+$ to the exact sequence
$$
0 \lra M(-1) \stackrel{x}{\lra} M \lra M/xM \lra 0
$$
shows that $f_{R_+}(M/xM) \geq r-1$ and
leads to an exact sequence of $R_0$-modules
$$
0 \lra H^{r-1}_{R_+}(M/xM)_n \lra H^r_{R_+}(M)_{n-1}
\stackrel{x}{\lra} H^r_{R_+}(M)_n
$$
for each $n \leq \widetilde{n}$.
The left exactness of the functor $\Gamma_{\fm_0}$ therefore leads to the
conclusion that, for each $x \in U$, the map
$$
J_{n-1} = \Gamma_{\fm_0}(H^r_{R_+}(M)_{n-1}) \stackrel{x}{\lra}
J_n = \Gamma_{\fm_0}(H^r_{R_+}(M)_n)
$$
is injective for all $n \leq \min \{\widetilde{n}, n_x\}$.  Hence $(0:_Jx)$
is an $R$-module of finite length, for all $x \in U$. Since $R_0/\fm_0$ is
algebraically closed, we can now deduce from \cite[Corollary (2.2)]{Brodm82}
that $J$ is an $R$-module of finite length, and this is a contradiction. We
have therefore proved that there exists $y \in U$ such that
$\fm_0 \in \Ass_{R_0}(H^{r-1}_{R_+}(M/yM)_n)$ for infinitely many $n < 0$.
This implies that $f_{R_+}(M/yM) \leq r-1$; therefore, as we have already
noted that $f_{R_+}(M/yM) \geq r-1$, we must have $f_{R_+}(M/yM) = r-1$.
Hence, by \cite[Proposition (5.6)]{BroHel00},
$\Ass_{R_0}(H^{r-1}_{R_+}(M/yM)_n)$ is
asymptotically stable for $n \rightarrow -\infty$; therefore
$\fm_0 \in \Ass_{R_0}(H^{r-1}_{R_+}(M/yM)_n)$ for all $n < < 0$.
\end{proof}

\begin{lem}\label{bo.44} The notation is as in {\rm \S \ref{in}} and
{\rm \ref{mi.0}}. In addition, suppose that $(R_0,\fm_0)$ is a
regular local ring of dimension $d$ and that $R = R_0[X_1, \ldots,
X_r]$, a polynomial ring graded in the usual way.

Assume that $f_{R_+}(M) < r$ and that $0 \rightarrow N \rightarrow F
\rightarrow M \rightarrow 0$ is an exact sequence of
finitely generated graded $R$-modules
and homogeneous homomorphisms in which $F$ is free. Then

\begin{enumerate}
\item $\depth N_{\fp} = \min\left\{\height \fp, \depth M_{\fp} + 1\right\}$
for all $\fp \in \Supp (N)$;
\item for $i \in \nn$,
the (necessarily homogeneous) connecting homomorphism $H^i_{R_+}(M)
\rightarrow H^{i+1}_{R_+}(N)$ induced by the given exact sequence is
an isomorphism when $i < r-1$ and a monomorphism when $i = r-1$; and
\item $f_{R_+}(N) = f_{R_+}(M) + 1$.
\end{enumerate}
\end{lem}

\begin{proof} Note that $N \neq 0$ because $f_{R_+}(F) = r$.

(i) This is immediate from the exact sequence
$0 \rightarrow N_{\fp} \rightarrow F_{\fp} \rightarrow M_{\fp} \rightarrow 0$.

(ii) This is immediate from the fact that $H^i_{R_+}(F) = 0$ for all
$i < r$.

(iii) This now follows from part (ii) and the hypothesis that
$f_{R_+}(M) < r$.
\end{proof}

\begin{lem}\label{bo.45} Assume that
$(R_0,\fm_)$ is a regular local ring. Then there exists a regular
local flat extension ring $(R'_0,\fm'_0)$ of $R_0$ such that
$\fm_0R'_0 = \fm'_0$ and $R'_0/\fm'_0$ is algebraically closed.
\end{lem}

\begin{proof} Denote as usual $\dim R_0$ by $d$. Let
$(\widehat{R_0}, \widehat{\fm_0})$ denote the
completion of $R_0$, so that 
$\widehat{\fm_0} = \fm_0\widehat{R_0}$; 
of course, this is a regular local flat extension ring
of $R_0$ of dimension $d$. By \cite[Proposition (2.2)]{BrMaMi00}, there exists
a (Noetherian) local
flat extension ring $(R'_0,\fm'_0)$ of $\widehat{R_0}$
such that $\widehat{\fm_0}R'_0 = \fm'_0$
and $R'_0/\fm'_0$ is algebraically closed. Therefore $\fm_0R'_0 = \fm'_0$,
so that $\fm'_0$ can be generated by $d$ elements. By flatness, $\dim R'_0 \geq
d$, and so $(R'_0,\fm'_0)$ is a regular local ring of dimension $d$.
\end{proof}

We are now ready to present our main result of this section.

\begin{thm}\label{bo.35} Assume that the graded ring $R$ is a
homomorphic image of a regular (commutative Noetherian) ring, and
that the non-zero graded $R$-module $M = \bigoplus_{n \in \Z} M_n$
is finitely generated and not $R_+$-torsion.  Set
$$
f := f_{R_+}(M) =
\inf \left\{ i \in \N : H^i_{R_+}(M)
\mbox{~is not finitely generated} \right\}.
$$
Then
$$
\Ass_{R_0}(H^f_{R_+}(M)_n) =
\left\{ \fp \cap R_0 : \fp \in \Proj (R) \text{~and~} \depth
M_{\fp} + \height (\fp + R_+)/\fp = f\right\} 
\quad \text{~for all~} n < < 0.
$$
\end{thm}

\begin{note} {\rm By Grothendieck's Finiteness Theorem (see
\cite[13.1.17]{BS}), the set on the right-hand side of the final
display in the statement of the theorem is non-empty:
note that $f$ is finite. A consequence of this theorem is that
that set is finite.}
\end{note}

\begin{proof} We first show
by induction on $f$ that, for $\fp \in \Proj (R)$
with $\height (\fp + R_+)/\fp= f$, we have $\fp \cap R_0 \in
\Ass_{R_0}(H^f_{R_+}(M)_n)$ for all $n < < 0$. Now $\height (\fp +
R_+)/\fp \geq 1$; so, if $f = 1$ and $\depth M_{\fp} + \height
(\fp + R_+)/\fp = 1$, then $\fp \in \Ass_RM$. The claim in the
case when $f = 1$ is therefore immediate from Lemma \ref{bo.14}.

Thus we assume now that $f > 1$ and make the obvious inductive
assumption.  One can use homogeneous localization at $\fp + R_+$
to see that it is enough to complete the inductive step under the
additional hypotheses that $R$ is $\mbox{\rm *}$local with unique
$\mbox{\rm *}$maximal ideal $\fm$, and that $\fm_0 := \fm \cap R_0
= \fp_0$.

Set $\overline{M} := M/\Gamma_{R_+}(M)$; recall (\cite[2.1.7]{BS})
that there are homogeneous isomorphisms $H^i_{R_+}(M)
\stackrel{\cong}{\lra} H^i_{R_+}(\overline{M})$ for each $i \in
\N$. Since $M_{\fp} \cong \overline{M}_{\fp}$, it follows that one
may assume, in this inductive step, that $\Gamma_{R_+}(M) = 0$.

The argument now splits into two cases, according as $\fp \in
\Ass_RM$ or $\fp \not\in \Ass_RM$. In the first case, it follows
from \cite[15.1.2]{BS} that there exists a positive integer $d$
and a homogeneous element $g_d \in R_d$ which is a non-zerodivisor
on $M$.  Let $\fq$ be a minimal prime ideal of $\fp + Rg_d$;
necessarily, $\fq \in \Proj (R)$ (since $f > 1$), and $\fq \cap
R_0 = \fm_0$. The catenarity of $R$ ensures that $\height
\left((\fq + R_+)/\fq\right) = f-1$. It follows from \cite[Chapter
6, Lemma 4]{HMold} that $\fq \in \Ass(M/g_dM)$, and so one can use
Grothendieck's Finiteness Theorem (see \cite[9.5.2]{BS}) to see
that $f_{R_+}(M/g_dM) \leq 0 + \height \left((\fq +
R_+)/\fq\right) = f-1$.

In the second case, when $\fp \not\in \Ass_RM$, we choose $g_d$ as
follows. First note that, for each $\fq' \in \Ass_RM$, we have $\fp
\cap R_+ \not\subseteq \fq'$. To see this, suppose that $\fp \cap
R_+ \subseteq \fq'$ for some $\fq' \in \Ass_RM$. Then $\fp \subset
\fq'$ (since $\Gamma_{R_+}(M) = 0$), so that (since $\fq' \cap R_0
\supseteq \fp\cap R_0 = \fm_0$), we have
$$
\height \left((\fq' + R_+)/\fq'\right) = \height
\left(\fm/\fq'\right)< \height \left(\fm/\fp\right) = \height
\left((\fp + R_+)/\fp\right).
$$
This implies that $\depth M_{\fq'} + \height \left((\fq' +
R_+)/\fq'\right) < f-1$, contrary to Grothendieck's Finiteness
Theorem.  We have therefore shown that $\fp \cap R_+ \not\subseteq
\fq'$. As this is true for all $\fq' \in \Ass_RM$, we can now use
\cite[15.1.2]{BS} to see that there exists a positive integer $d$
and a homogeneous element $g_d \in \fp \cap R_d$ which is a
non-zerodivisor on $M$. Note that $\depth (M/g_dM)_{\fp} = \depth
M_{\fp} - 1$, so that, by Grothendieck's Finiteness Theorem,
$$
f_{R_+}(M/g_dM) \leq \depth (M/g_dM)_{\fp} + \height \left((\fp +
R_+)/\fp\right) = f-1.
$$

Thus, in both cases, we have found a homogeneous element $g_d$ of
$R$ of positive degree $d$ which is a non-zerodivsor on $M$ and is
such that $f_{R_+}(M/g_dM) \leq f-1$. Application of local
cohomology with respect to $R_+$ to the exact sequence
$$
0 \lra M(-d) \stackrel{g_d}{\lra} M \lra M/g_dM \lra 0
$$
shows that $f_{R_+}(M/g_dM) \geq f-1$, and that, for all $n < <
0$, the $R_0$-module $H^f_{R_+}(M)_n$ has a submodule isomorphic
to $H^{f-1}_{R_+}(M/g_dM)_{n+d}$. It therefore follows that
$f_{R_+}(M/g_dM) = f-1$ (so that $M/g_dM$ is not $R_+$-torsion),
and we can apply the inductive hypothesis to $M/g_dM$.

In our first case, when $\fp \in \Ass_RM$, we have already noted
that $\fq \in \Proj (R) \cap \Ass(M/g_dM)$, that $\fq \cap R_0 =
\fm_0$, and that $\depth (M/g_dM)_{\fq} + \height \left((\fq +
R_+)/\fq\right) = f-1$. We therefore use $\fq$ to draw a
conclusion from the inductive hypothesis.

In our second case, when $\fp \not\in \Ass_RM$, we noted that
$\depth (M/g_dM)_{\fp} + \height \left((\fp + R_+)/\fp\right) =
f-1$; in this case, we use $\fp$ to draw a conclusion from the
inductive hypothesis.

In both cases, the inductive hypothesis yields that $\fm_0 \in
\Ass_{R_0}(H^{f-1}_{R_+}(M/g_dM)_{n+d})$ for all $n < < 0$.
Therefore $\fm_0 \in \Ass_{R_0}(H^{f}_{R_+}(M)_{n})$ for all $n <
< 0$, and the inductive step is complete.

We have thus proved that
$$
\Ass_{R_0}(H^f_{R_+}(M)_n) \supseteq
\left\{ \fp \cap R_0 : \fp \in \Proj (R) \text{~and~} \depth
M_{\fp} + \height (\fp + R_+)/\fp = f\right\} 
\quad \text{~for all~} n < < 0.
$$

To complete the proof, we suppose that $\fp_0 \in
\Ass_{R_0}(H^f_{R_+}(M)_n)$ for all $n < < 0$; it is enough for us to show that
there exists
$\fp \in \Proj (R)$ with $\fp \cap R_0 = \fp_0$ and $\depth
M_{\fp} + \height (\fp + R_+)/\fp = f$. Our first steps in this direction
show that additional simplifications are posssible.

Invert $R_0 \setminus \fp_0$; in other words, apply homogeneous
localization at $\fp_0 + R_+$. Observe that the hypotheses imply
that $(R_0)_{\fp_0}$ is a homomorphic image of a regular local
ring $(R'_0,\fm'_0)$ and that $R_{(\fp_0 + R_+)}$ is an image of a
polynomial ring $R' := R'_0[X_1, \ldots, X_r]$, graded in the
usual way, under a ring homomorphism which is homogeneous in the
sense of \cite[Definition 13.1.2]{BS}. Consider $M$ as a finitely
generated graded $R'$-module; we can then use the Graded
Independence Theorem \cite[13.1.6]{BS} to see that $f =
f_{R'_+}(M)$ and that it is enough for us establish the existence
of a $\fp \in \Proj(R)$ with the specified properties under the
additional hypotheses that $(R_0,\fm_0)$ is a regular local ring,
that $\fp_0 = \fm_0$, and that $R = R_0[X_1, \ldots, X_r]$, a
polynomial ring graded in the usual way.

We deal first with the case where $r = 1$. Then $f = 1$.
Set $\overline{M} := M/\Gamma_{R_+}(M)$, and recall that there is
a homogeneous isomorphism
$H^1_{R_+}(M) \stackrel{\cong}{\lra} H^1_{R_+}(\overline{M})$.
Therefore, since $M_{\fq} \cong \overline{M}_{\fq}$ for all $\fq
\in \Spec (R) \setminus \Var (R_+)$, we can, and do, impose
the additional hypothesis that $\Gamma_{R_+}(M) = 0$. (Here, $\Var (\fa)$, 
for an
ideal $\fa$ of $R$, denotes the variety of $\fa$.)

By hypothesis, $\left( 0 :_{H^1_{R_+}(M)} :\fm_0 \right)_n \neq 0$ for all
$n < < 0$. It therefore follows that the graded $R$-module
$\Gamma_{\fm_0R}\left(H^1_{R_+}(M)\right)$ is not finitely generated.

Let
\[
0
\longrightarrow \mbox{\rm *}E^0(M)
\stackrel{d^0}{\longrightarrow} \mbox{\rm *}E^1(M)
\stackrel{d^1}{\longrightarrow} \mbox{\rm *}E^{2}(M)
\longrightarrow \cdots \longrightarrow
\mbox{\rm *}E^i(M) \longrightarrow \cdots
\]
be the minimal $\mbox{\rm *}$injective resolution of $M$, with associated
(necessarily homogeneous) augmentation homomorphism $d^{-1} : M \longrightarrow
\mbox{\rm *}E^0(M)$. Since $\Gamma_{R_+}(M) = 0$, it follows from
\cite[Theorem 2.4]{70} that $\Gamma_{R_+}(\mbox{\rm *}E^0(M)) = 0$,
so that $\Gamma_{\fm}(\mbox{\rm *}E^0(M)) = 0$. Therefore
$$
H^1_{R_+}(M) \cong \Ker \left(\Gamma_{R_+}(d^1)\right) \quad \mbox{~and~} \quad
H^1_{\fm}(M) \cong \Ker \left(\Gamma_{\fm}(d^1)\right).
$$
Here, $\Gamma_{R_+}(d^1) : \Gamma_{R_+}\!\left(\mbox{\rm
*}E^1(M)\right) \lra \Gamma_{R_+}\!\left(\mbox{\rm *}E^{2}(M)\right)$
is the map induced by $d^1$, {\it etcetera}. Thus
$$
\Gamma_{\fm_0R}\left(H^1_{R_+}(M)\right) \cong
\Gamma_{\fm_0R}\left(\Ker \left(\Gamma_{R_+}(d^1)\right)\right) =
\Ker \left(\Gamma_{\fm}(d^1)\right) \cong H^1_{\fm}(M).
$$
Therefore, $H^1_{\fm}(M)$ is not finitely
generated. Hence, by Grothendieck's Finiteness Theorem (see
\cite[13.1.17]{BS}), there exists $\fp \in
\mbox{\rm *}\Spec (R) \setminus \Var
(\fm)$ such that $\depth M_{\fp} + \height \fm/\fp = 1$.
This means that
$\fp \in \Ass_RM$ and $\height \fm/\fp = 1$.  Note that
$\fp \not\supseteq R_+$, because $\Gamma_{R_+}(M) = 0$. Therefore
$\fp_0 := \fp \cap R_0 = \fm_0$, since otherwise $\fm \supset
\fp_0 + R_+ \supset \fp$ would be a chain of distinct prime ideals
of $R$, contrary to the fact that $\height \fm/\fp =
1$. The claim is therefore proved in the case where $r = 1$.

Now suppose that $r \geq 2$, and that the desired result has been proved
for smaller values of $r$. Note that, by Grothendieck's Finiteness Theorem,
it is enough for us to show that there exists
$\fp \in \Proj (R)\cap \Var(\fm_0R)$ with $\depth
M_{\fp} + \height (\fp + R_+)/\fp \leq f$.

By Lemma \ref{bo.45}, there exists a regular local flat extension
ring $(R'_0,\fm'_0)$ of $R_0$ such that $\fm_0R'_0 = \fm'_0$ and
$R'_0/\fm'_0$ is algebraically closed. Let $R' = R
\otimes_{R_0}R'_0$, which we identify with $R'_0[X_1, \ldots,
X_r]$ in the obvious way. Let $M'$ denote the finitely generated
graded $R'$-module $M\otimes_RR'$. It follows from \cite[13.1.8
and 15.2.2]{BS} that $f_{R'_+}(M') = f$ and $\fm'_0 \in
\Ass_{R'_0}(H^f_{R'_+}(M')_n)$ for all $n < < 0$.

Suppose that we have found
$\fp' \in \Proj (R')\cap \Var (\fm'_0R')$ such that
$\depth M'_{\fp'} + \height (\fp' + R'_+)/\fp' \leq f$.
Set $\fp := \fp' \cap R$. Then it follows from Lemma \ref{bo.42} that
$\fp \in \Proj (R)\cap \Var(\fm_0R)$ and
$$
\depth M_{\fp} + \height (\fp + R_+)/\fp \leq
\depth_{R'_{\fp'}} M'_{\fp'} + \height (\fp' + R'_+)/\fp' \leq f.
$$
Therefore we can, and do, assume for the remainder of this proof that
$R_0/\fm_0$ is algebraically closed.

We now proceed by descending induction on $f$.  Note that $f \leq r$. So
we deal first with the case where $f = r$. By Lemma \ref{bo.43},
there exists $y_r \in
R_1 \setminus \fm_0R_1$ such that $y_r$
is a non-zerodivisor on $M/\Gamma_{R_+}(M)$,
that $f_{R_+}(M/y_rM) = r-1$ and that
$$
\fm_0 \in
\Ass_{R_0}(H^{r-1}_{R_+}(M/y_rM)_n) \quad \text{~for all~} n < < 0.
$$
Since the image of $y_r$ in $R_1/\fm_0R_1$ is non-zero, there exist
$y_1, \ldots, y_{r-1} \in R_1$ such that $R_1$ is generated (over $R_0$) by
$y_1, \ldots, y_{r-1},y_r$. Note that $R = R_0[y_1, \ldots, y_{r-1},y_r]$ and
that $y_1, \ldots, y_{r-1},y_r$ are algebraically independent over $R_0$.
Therefore, we can, and do, assume that $y_r = X_r$.

We can consider $M/X_rM$ as a finitely generated graded module over $R/X_rR$.
The Graded Independence Theorem \cite[13.1.6]{BS} shows that
$f_{(R/X_rR)_+}(M/X_rM) = r-1$ and that
$$
\fm_0 \in
\Ass_{R_0}(H^{r-1}_{(R/X_rR)_+}(M/X_rM)_n) \quad \text{~for all~} n < < 0.
$$
Since $R/X_rR$ is (homogeneously) isomorphic to $R_0[X_1
,\ldots,X_{r-1}]$, we may apply the inductive hypothesis to deduce
that there exists $\overline{\fp} \in \Proj (R/X_rR)\cap
\Var(\fm_0(R/X_rR))$ with $$\depth (M/X_rM)_{\overline{\fp}} +
\height (\overline{\fp} + (R/X_rR)_+)/\overline{\fp} \leq r-1.$$
Let $\fp$ be the inverse image of $\overline{\fp}$ under the
natural ring homomorphism $R \rightarrow R/X_rR$. Then $X_r \in
\fp$ and $\fp \in \Proj (R)\cap \Var(\fm_0R)$; also $\depth
M_{\fp} = \depth (M/X_rM)_{\overline{\fp}} + 1$ (because
$M_{\fp} \cong (M/\Gamma_{R_+}(M))_{\fp}$) and
$$
\height (\fp + R_+)/\fp = \height (\overline{\fp} + (R/X_rR)_+)/\overline{\fp}.
$$
Hence $\depth M_{\fp} + \height (\fp + R_+)/\fp \leq r$.
Thus we have found a $\fp$ with the required properties in the case
where $f = r > 1$.

Now suppose that $f < r$ and that the desired result has been proved for
larger values of $f$ (for this value of $r$). There is an exact sequence
$0 \rightarrow N \rightarrow F
\rightarrow M \rightarrow 0$ of
finitely generated graded $R$-modules
and homogeneous homomorphisms in which $F$ is free. By Lemma \ref{bo.44}(iii),
we have $f_{R_+}(N) = f + 1$; by part (ii) of the same lemma,
$\fm_0 \in \Ass_{R_0}(H^{f+1}_{R_+}(N)_n)$ for all $n < < 0$. Therefore, by
the inductive hypothesis,
there exists
$\fp \in \Proj (R)\cap \Var(\fm_0R)$ with
$\depth N_{\fp} + \height (\fp + R_+)/\fp \leq f+1$.

Note that $\depth M_{\fp} \leq \height \fp$: we consider the cases where
$\depth M_{\fp} < \height \fp$ and $\depth M_{\fp} = \height \fp$
separately.
When $\depth M_{\fp} < \height \fp$, it follows from Lemma \ref{bo.44}(i)
that $\depth N_{\fp} = \depth M_{\fp} + 1$; therefore
$\depth M_{\fp} + \height (\fp + R_+)/\fp \leq f$. In the other case,
$\depth M_{\fp} = \height \fp$, so that, again by Lemma \ref{bo.44}(i),
$\depth N_{\fp} = \height \fp$. Therefore, in this case,
\begin{align*}
\height (\fm_0 + R_+) & = \height (\fp + R_+) \\
& = \height (\fp + R_+)/\fp + \height \fp \\
& = \height (\fp + R_+)/\fp + \depth N_{\fp} \leq f+1 \leq r.
\end{align*}
Therefore $\fm_0 = 0$ and $R_0$ is a field. In this case, the
desired conclusion is clear from the graded version of Grothendieck's
Finiteness Theorem (see \cite[13.1.17]{BS}).
The proof is now complete.
\end{proof}

\section{\bf Further examination of Singh's example}
\label{se}

In \cite[\S 4]{Singh00}, A. K. Singh showed that the ring $R' :=
\Z[X,Y,Z,U,V,W]/(XU + YV + ZW)$, where $X,Y,Z,U,V,W$ are
independent indeterminates over $\Z$, has the property that
$\Ass_{R'}(H^3_{\fa}(R'))$ is infinite, where $\fa$ is the ideal
generated by the images of $U, V, W$. Brodmann and Hellus
\cite[(5.7)(A)]{BroHel00} observed that Singh's argument leads to
an interesting conclusion about graded components of graded local
cohomology modules: we can consider $\Z[X,Y,Z,U,V,W]$ as a
positively graded ring with $0$th component $\Z[X,Y,Z]$ and
$U,V,W$ each assigned degree $1$; $R'$ inherits a structure as a
standard positively graded ring with
$R'_+ = \fa$; the argument Singh used to prove his
result mentioned above actually shows that $\left\{ \fp \cap \Z :
\fp \in \Ass_{R'}(H^3_{R'_+}(R'))\right\}$ is an infinite set, and
Brodmann and Hellus noted that this implies that
$\Ass_{R'_0}(H^3_{R'_+}(R')_n)$ is not asymptotically
stable for $n \lra \infty$.

Our aim in the rest of this paper is to use Gr\"obner basis
techniques on Singh's example to identify precisely the set
$\Ass_{R'_0}(H^3_{R'_+}(R')_n)$ for each $n \leq -3$, and to then
deduce that $\Ass_{R'_0}(H^3_{R'_+}(R')_n)$ is not asymptotically
increasing for $n \lra \infty$.

\begin{ntn}\label{se.1} {\rm Throughout the rest of the paper, the symbol $L$
will denote either a field or a principal ideal domain (PID), and
$R$ will denote the polynomial ring $L[X,Y,Z,U,V,W]$, graded so
that $U,V,W$ have degree $1$ and $X,Y,Z$ have degree $0$; thus
$R_0 = L[X,Y,Z]$. We shall set $F := XU + YV + ZW$, and $R' :=
R/FR$, again a standard positively graded ring. The natural map $R
\rightarrow R'$ maps $R_0$ isomorphically onto $R'_0$, and so we
shall identify elements of $R_0$ with their natural images in
$R'_0$.  In the case where $L = \Z$, the rings $R$ and $R'$ are
those occurring in Singh's example mentioned above. However, it
will be helpful in another context to have some calculations
available in the case where $L$ is the rational field, for
example.

Since $H^3_{R'_+}(R')$ is homogeneously isomorphic to
$H^3_{R_+}(R/FR)$, we can use the exact sequence
$$
H^3_{R_+}(R)(-1) \stackrel{F}{\lra} H^3_{R_+}(R) \lra
H^3_{R_+}(R/FR) \lra 0
$$
of graded $R$-modules and homogeneous homomorphisms (induced from
the exact sequence
$$
0 \lra R(-1) \stackrel{F}{\lra} R \lra R/FR \lra 0 )
$$
to study $H^3_{R'_+}(R')$. Furthermore, we can realize
$H^3_{R_+}(R)$ as the module $R_0[U^-,V^-,W^-]$ of inverse
polynomials described in \cite[12.4.1]{BS}: this graded $R$-module
has end $-3$, and, for each $d \geq 3$, its $(-d)$-th component
is a free $R_0$-module of rank ${{d-1} \choose {2}}$ with base
$
\left(U^{\alpha}V^{\beta}W^{\gamma}\right)_{-\alpha,
-\beta,-\gamma \in \N,~\alpha + \beta + \gamma = -d}.
$
We plan to study the graded components of $H^3_{R_+}(R/FR)$ by considering
the cokernels of the $R_0$-homomorphisms
$$
F_{-d}: R_0[U^-,V^-,W^-]_{-d-1} \lra R_0[U^-,V^-,W^-]_{-d} \quad (d \geq 3)
$$
given by multiplication by $F$. In order to represent these
$R_0$-homomorphisms between free $R_0$-modules by matrices, we
specify an ordering for each of the above-mentioned bases by
declaring that $$U^{\alpha_1}V^{\beta_1}W^{\gamma_1} <
U^{\alpha_2}V^{\beta_2}W^{\gamma_2}$$ (where $-\alpha_i,
-\beta_i,-\gamma_i \in \N$ and $\alpha_i + \beta_i + \gamma_i = n
\leq -3$ for $i = 1,2$) precisely when $\alpha_1 > \alpha_2$ or
$\alpha_1 = \alpha_2$ and $\beta_1 > \beta_2$. For example, this
ordering on our base for $R_0[U^-,V^-,W^-]_{-5}$ is such that
$$
U^{-1}V^{-1}W^{-3} < U^{-1}V^{-2}W^{-2} < U^{-1}V^{-3}W^{-1} <
U^{-2}V^{-1}W^{-2} < U^{-2}V^{-2}W^{-1} < U^{-3}V^{-1}W^{-1}.
$$

We shall frequently need to consider an $R_0$-homomorphism from
the free $R_0$-module $R_0^n$ (regarded as consisting of column
vectors) to $R_0^m$ (where $m$ and $n$ are positive integers)
given by left multiplication by an $m \times n$ matrix $C$ with
entries in $R_0$. In these circumstances, we shall also use $C$ to
denote the homomorphism; its image $\Ima C$ is just the submodule
of $R_0^m$ generated by the columns of $C$, for if
$(\mathbf{e}_i)_{i=1,\ldots,n}$ denotes the standard base for
$R_0^n$, then $C\mathbf{e}_j$ is just the $j$-th column of $C$
(for $1 \leq j \leq n$).

The theory of Gr\"obner bases is well developed for ideals in
polynomial rings in finitely many indeterminates with coefficients
in a principal ideal domain, and for submodules of finite free
modules over polynomial rings in finitely many indeterminates over
a field (see, for example, \cite[Chapter 3]{AL}). It is
straightforward to combine the methods from these two parts of the
theory to produce a theory of Gr\"obner bases for submodules of
finite free $L[X,Y,Z]$-modules. Thus much of the work below
applies both to the case where $L$ is a PID and the case where $L$
is a field.

In this paper, we use the lexicographical term order with $X > Y > Z$
in $R_0$, and for each $n \in \N$ we set $>$ to be the `term-over-position'
extension of this order to $R_0^n$ defined as follows: a monomial in
$R_0^n$ is a column vector of the form $m\mathbf{e}_j$, where $m$ is a
monomial in $R_0$ and $\mathbf{e}_j$ is the $j$-th standard base vector of
$R_0^n$; and $m_1\mathbf{e}_{j_1} > m_2\mathbf{e}_{j_2}$ (for monomials $m_1,
m_2$ of $R_0$ and $j_1,j_2 \in \{1, \ldots, n\}$) if and only if
$$
m_1 > m_2 \quad \mbox{or} \quad m_1 = m_2 \mbox{~and~} j_1 < j_2.
$$

If $A$ is an $m \times n$ matrix with entries in $R_0$ (we shall
say `over $R_0$') and $\mathbf{f}, \mathbf{h} \in R_0^n$, then we
shall say that \textit{$\mathbf{f}$ reduces to $\mathbf{h}$ modulo
$A$}, denoted by $\mathbf{f} \stackrel{A}{\lra}_+ \mathbf{h}$,
when $\mathbf{f}$ reduces to $\mathbf{h}$ modulo the set of
columns of $A$ (see \cite[Definition 3.5.8]{AL}, but modify that
definition to imitate \cite[Definitions 4.1.1 and 4.1.6]{AL} in
the case where $L$ is a PID).

We shall denote the leading monomial, leading coefficient and
leading term of $\mathbf{f} \in R_0^n$ by $\lm(\mathbf{f})$,
$\lc(\mathbf{f})$ and $\lt(\mathbf{f})$ respectively.

We shall use $I_n$ to denote the $n \times n$ identity matrix. For
each $n \in \N$, we let $A_n$ denote the $n \times (n+1)$ matrix
given by
$$
A_n = \left[\begin{array}{ccccc}
  Z & Y & 0 & \dots & 0 \\
  0 & Z & Y & 0 & \dots \\
  \vdots &  & \ddots & \ddots &  \\
  0 & \dots & 0 & Z & Y
\end{array}\right].
$$ }
\end{ntn}

\begin{lem}\label{se.1a} Let $d \in \N$ with $d \geq 3$.

\begin{enumerate}
\item With the
notation of\/ {\rm \ref{se.1}}, the
$R_0$-homomorphism
$$
F_{-d}: R_0[U^-,V^-,W^-]_{-d-1} \lra R_0[U^-,V^-,W^-]_{-d}
$$
given by multiplication by $F$ is represented, relative to the
bases specified in\/ {\rm \ref{se.1}} listed in increasing order,
by the ${{d-1} \choose {2}} \times {{d} \choose {2}}$  matrix
$$T_d :=
\left[ \begin{array}{cccccc}
  A_{d-2} & XI_{d-2} & 0 & \dots &  & 0 \\
  0 & A_{d-3} & XI_{d-3} & 0 & \dots & 0 \\
  0 & 0 & A_{d-4} & XI_{d-4} &  & 0 \\
  \vdots & \vdots &  & \ddots & \ddots & \vdots \\
  0 & 0 & \dots & 0 & A_1 & XI_1
\end{array} \right],
$$
where $A_{d-2}, \ldots, A_1$ are as defined in\/ {\rm \ref{se.1}}.

\item Each associated prime in
$\Ass_{R'_0}(H^3_{R'_+}(R')_{-d})$ contains $X$, $Y$ and $Z$.

\item We have $(X,Y,Z) \in \Ass_{R'_0}(H^3_{R'_+}(R')_{-d})$.
\end{enumerate}
\end{lem}

\begin{proof} (i) This
follows from the fact that, for negative integers
$\alpha,\beta,\gamma$,
$$
F (U^\alpha V^\beta W^\gamma) =
X(1-\delta_{\alpha, -1}) U^{\alpha+1} V^\beta W^\gamma+
Y(1-\delta_{\beta, -1}) U^{\alpha} V^{\beta+1} W^\gamma+
Z(1-\delta_{\gamma, -1}) U^{\alpha} V^\beta W^{\gamma+1},
$$
where $\delta_{i,j}$ is Kronecker's delta.

(ii) Consider the last column of $T_d$ to see that
$X\mathbf{e}_{{d-1 \choose 2}} \in \Ima T_d$; therefore
$XY\mathbf{e}_{{d-1 \choose 2}}, XZ\mathbf{e}_{{d-1 \choose 2}}
\in \Ima T_d$, so that $X^2\mathbf{e}_{{d-1 \choose 2}-1},
X^2\mathbf{e}_{{d-1 \choose 2}-2} \in \Ima T_d$ in view of the
next-to-last and second-to-last columns of $T_d$; we can now
continue in this way to see that each element of $\Coker T_d =
\Coker F_{-d}$ is annihilated by $X^{d-2}$.  By symmetry,
$Y^{d-2}$ and $Z^{d-2}$ also annihilate $\Coker F_{-d} = \Coker
T_d$.

(iii) It is clear from part (i) that $(\Ima F_{-3}:_{R_0}
U^{-1}V^{-1}W^{-1}) = (X,Y,Z)$. Hence $(0 :_{R_0} \Coker F_{-3}) =
(X,Y,Z)$. Now multiplication by $U^{d-3}$ induces an
$R_0$-epimorphism $ \Coker F_{-d} \lra \Coker F_{-3} $, so that,
in view of the above proof of part (ii), we have
$$
(X^{d-2},Y^{d-2},Z^{d-2}) \subseteq (0 :_{R_0} \Coker F_{-d}) \subseteq
(0 :_{R_0} \Coker F_{-3}) = (X,Y,Z).
$$
Therefore $(X,Y,Z)$ is a minimal member of the support of $\Coker
F_{-d}$.
\end{proof}

\begin{lem}\label{se.2} Let $k,m,n,q \in \nn$ with $m,n > 0$. Let $A = \left[
a_{ij}\right]$ be an $m \times n$ matrix with entries in $L[Y,Z]$,
let $\mathbf{f} \in L[Y,Z]^n$, and let $M$ and $M'$ denote the $(k
+ n + m + q)$-rowed block matrices over $R_0$ given by
$$
M := \left[ \begin{array}{cc} 0 & 0 \\ XI_n & \mathbf{f}
\\ A & 0 \\ 0 & 0 \end{array}
\right] \quad \mbox{~and~} \quad
M' := \left[ \begin{array}{c} 0 \\ XI_n
\\ A \\ 0 \end{array}
\right],
$$
in which the first $k$ and last $q$ rows are all zero.

Then each $S$-polynomial of two columns of $M$ is either $0$ or reduces
modulo $M'$ to
$$\left[ \begin{array}{c} 0 \\ 0 \\ \pm A\mathbf{f}\\ 0 \end{array} \right]$$
(in which the lowest `$0$' stands for the $q \times 1$ zero
matrix), and these column matrices do arise from $S$-polynomials
in this way.
\end{lem}

\begin{proof} Suppose $\mathbf{f} \neq 0$, and let
$\mathbf{f} = \sum_{j=1} ^{t} c_{i_j}T_{i_j}\mathbf{e}_{i_j}$ be
an expression for $\mathbf{f}$ as a sum of terms, where the
$T_{i_j}$ are monomials in $Y$ and $Z$ and the $c_{i_j}$ are
elements of $L$; suppose that $\lt(\mathbf{f}) =
c_{i_h}T_{i_h}\mathbf{e}_{i_h}$. Let $\mathbf{m}_{j}$ denote the
$j$th column of $M$, for each $j = 1, \ldots, n+1$.

Since $\mathbf{f} \in L[Y,Z]$, we have $\lcm (T_{i_h},X) =
T_{i_h}X$. All $S$-polynomials of two columns of $M$ are zero
except possibly for those of $\mathbf{m}_{i_h}$ and
$\mathbf{m}_{n+1}$.  Note that
$$
\mathbf{m}_{n+1} = \sum_{j=1}^{t} c_{i_j}T_{i_j}\mathbf{e}_{i_j+k}
\quad \mbox{~and~} \quad \mathbf{m}_{i} = X\mathbf{e}_{i+k} +
\sum_{\rho=1}^{m}a_{\rho i}\mathbf{e}_{\rho+k+n} ~(1 \leq i \leq
n).
$$
We have
\begin{align*}
S(\mathbf{m}_{i_h},\mathbf{m}_{n+1}) & =
\frac{c_{i_h}}{1}\frac{T_{i_h}X}{X}\mathbf{m}_{i_h}
- \frac{c_{i_h}}{c_{i_h}}\frac{T_{i_h}X}{T_{i_h}}\mathbf{m}_{n+1}
= c_{i_h}T_{i_h}\mathbf{m}_{i_h} - X\mathbf{m}_{n+1} \\
& =
c_{i_h}T_{i_h}X\mathbf{e}_{i_h+k} +
\sum_{\rho=1}^{m}a_{\rho i_h}c_{i_h}T_{i_h}\mathbf{e}_{\rho+k+n}
- \sum_{j=1}^{t} c_{i_j}XT_{i_j}\mathbf{e}_{i_j+k}\\
& =
\sum_{\rho = 1}^m a_{\rho i_h} c_{i_h}T_{i_h}\mathbf{e}_{\rho + k +n}
- \sum_{\stackrel{\scriptstyle{j=1}}{j \neq h}}^{t}
c_{i_j}XT_{i_j}\mathbf{e}_{i_j+ k}
\\
& \stackrel{M'}{\lra}_+ \sum_{\rho = 1}^m    a_{\rho i_h}
c_{i_h}T_{i_h}\mathbf{e}_{\rho + k +n} -
\sum_{\stackrel{\scriptstyle{j=1}}{j \neq h}}^{t}
c_{i_j}XT_{i_j}\mathbf{e}_{i_j+k} +
\sum_{\stackrel{\scriptstyle{j=1}}{j \neq h}}^{t}
c_{i_j}T_{i_j}\mathbf{m}_{i_j} \\
& =
\sum_{j=1}^{t}
\sum_{\rho = 1}^m a_{\rho i_j}c_{i_j}T_{i_j}\mathbf{e}_{\rho + k +n}
= \left[ \begin{array}{c} 0 \\ 0 \\ A\mathbf{f}\\ 0 \end{array} \right],
\end{align*}
as claimed.
\end{proof}

\begin{thm}\label{se.3} Consider the matrix
$$T_d :=
\left[ \begin{array}{cccccc}
  A_{d-2} & XI_{d-2} & 0 & \dots &  & 0 \\
  0 & A_{d-3} & XI_{d-3} & 0 & \dots & 0 \\
  0 & 0 & A_{d-4} & XI_{d-4} &  & 0 \\
  \vdots & \vdots &  & \ddots & \ddots & \vdots \\
  0 & 0 & \dots & 0 & A_1 & XI_1
\end{array} \right]
$$
of\/ {\rm \ref{se.1a}}. Define matrices $G_{d-2}, G_{d-3}, \ldots,
G_1$ by descending induction as follows: let $G_{d-2}$ be a
$(d-2)$-rowed matrix with entries in $L[Y,Z]$ whose columns
include those of $A_{d-2}$ and provide a Gr\"obner basis for $\Ima
A_{d-2}$; for $i \in \N$ with $d-2 > i \geq 1$, on the assumption
that $G_{i+1}$ has been defined as an $(i+1)$-rowed matrix with
entries in $L[Y,Z]$, let $G_i$ be an $i$-rowed matrix with entries
in $L[Y,Z]$ whose columns include those of $A_iG_{i+1}$ and
provide a Gr\"obner basis for $\Ima A_iG_{i+1}$. Then
\begin{enumerate}
\item the columns of
$$T_d' :=
\left[ \begin{array}{cccccc|ccccc}
  A_{d-2} & XI_{d-2} & 0 & \dots &  & 0 & G_{d-2} & 0 & \dots &  & 0\\
  0 & A_{d-3} & XI_{d-3} & 0 & \dots & 0 & 0 & G_{d-3} & 0 & \dots & 0\\
  0 & 0 & A_{d-4} & XI_{d-4} &  & 0 & 0 & 0 & G_{d-4} &   & 0 \\
  \vdots & \vdots & & \ddots & \ddots & \vdots & \vdots & \vdots
  & & \ddots & \vdots  \\
  0 & 0 & \dots & 0 & A_1 & XI_1 & 0 & 0 & \dots & 0 & G_1
\end{array} \right]
$$
form a Gr\"obner basis for $\Ima T_d' = \Ima T_d$; and
\item the columns of
$$ H_d :=
\left[\begin{array}{ccccc}
A_{d-2} & 0                & 0 &  \ldots         &  0  \\
0       & A_{d-3} A_{d-2}  & 0 &  \ldots         &  0\\
0 & 0& A_{d-4}A_{d-3} A_{d-2} & \ldots           & 0\\
\vdots  &   \vdots    &             & \ddots    & \vdots \\

0        &  0    &  \ldots           &           & A_{1} A_{2}
\dots A_{d-2}
\end{array}\right]$$
generate $\Ima T_d \cap L[Y,Z]^{{d-1} \choose {2}}$.
\end{enumerate}
\end{thm}

\begin{proof} (i) Let $\mathbf{s} = S(\mathbf{f},\mathbf{g})$
be a non-zero $S$-polynomial of two
columns $\mathbf{f}$ and $\mathbf{g}$ of $T_{d}'$. There are
various cases to consider.

First of all, if $\mathbf{f}$ and $\mathbf{g}$ have leading terms
in one of the first $d-2$ rows of $T_d'$, then either $\mathbf{s}
\stackrel{T_{d}'}{\lra}_+ 0$ because the columns of $\left[A_{d-2}
| G_{d-2}\right]$ form a Gr\"obner basis, or else $\mathbf{s}$
reduces modulo $T_d'$ to a column of
$$
\left[\begin{array}{c}  0 \\ \pm A_{d-3}G_{d-2} \\ 0 \\ \vdots \\
0
\end{array}\right]
$$
by Lemma \ref{se.2}; since the columns of $G_{d-3}$ include the
columns of $A_{d-3}G_{d-2}$, it follows that $\mathbf{s}
\stackrel{T_{d}'}{\lra}_+ 0$ in this case also.

Now suppose that $i \in \N$ with $d-2 > i > 1$ and that
$\mathbf{f}$ and $\mathbf{g}$ have leading terms in row $$k +
\sum_{j=i+1}^{d-2} j \quad \mbox{~for some~} k \in
\{1,\ldots,i\}.$$ Then either $\mathbf{s}
\stackrel{T_{d}'}{\lra}_+ 0$ because the columns of $G_{i}$ form a
Gr\"obner basis, or else $\mathbf{s}$ reduces modulo $T_d'$ to a
column of
$$
 \left[ \begin{array}{c}0 \\ \vdots \\ 0 \\ \pm A_{i-1}G_i\\0\\\vdots \\ 0
\end{array}\right]
$$
(where the block $A_{i-1}G_i$ is in the rows corresponding to
those where the blocks $A_{i-1}$ and $G_{i-1}$ are positioned in
$T_d'$), by Lemma \ref{se.2} again; since the columns of $G_{i-1}$
include the columns of $A_{i-1}G_{i}$, it follows that $\mathbf{s}
\stackrel{T_{d}'}{\lra}_+ 0$ in this case also.

Finally, suppose that $\mathbf{f}$ and $\mathbf{g}$ have leading
terms in the last row of $T_d'$. In this case, either $\mathbf{s}
\stackrel{T_{d}'}{\lra}_+ 0$ because the columns of $G_{1}$ form a
Gr\"obner basis, or $\mathbf{s} = S\left(X\mathbf{e}_{{{d-1}
\choose {2}}},h\mathbf{e}_{{{d-1} \choose {2}}}\right)$ for some
$h \in L[Y,Z]$; in the latter case, $\mathbf{s}$ reduces to $0$ in
one step modulo $\left\{X\mathbf{e}_{{{d-1} \choose
{2}}}\right\}$.

Thus, in all cases, $\mathbf{s} \stackrel{T_{d}'}{\lra}_+ 0$.
Hence, by (the analogue of) \cite[Theorem 3.5.19]{AL}, the columns
of $T_d'$ form a Gr\"obner basis.

To complete the proof of part (i), it only remains for us to show
that $\Ima T_d' = \Ima T_d$. It is easy to see by descending
induction on $i$ that, for each $i = d-2, d-3, \ldots, 1$, the
columns of $G_i$ include the columns of $A_iA_{i+1}\ldots A_{d-2}$
and form a Gr\"obner basis (over $\Z[Y,Z]$) for $\Ima
A_iA_{i+1}\ldots A_{d-2}$; hence (over both $\Z[Y,Z]$ and
$\Z[X,Y,Z]$)
$$
\Ima A_iA_{i+1}\ldots A_{d-2} = \Ima A_iG_{i+1} = \Ima G_i
\quad \mbox{~for all~} i = d-3, d-4, \ldots, 1.
$$
By Lemma \ref{se.2}, for such an $i$, each column of
$$
\left[ \begin{array}{c} 0 \\ \vdots \\ 0 \\ A_iG_{i+1} \\ 0 \\ \vdots \\ 0
\end{array} \right]
$$
(in which the block $A_iG_{i+1}$ occupies the rows corresponding to
those occupied by $A_i$ in $T_d$) can by obtained as a result of
reducing modulo $T_d$ the $S$-polynomial of a column of $T_d$ and a
column of
$$
\left[ \begin{array}{c} 0 \\ \vdots \\ 0 \\ G_{i+1} \\ 0 \\ \vdots \\ 0
\end{array} \right]
$$
(in which the block $G_{i+1}$ occupies the rows corresponding to
those occupied by $A_{i+1}$ in $T_d$). The claim in part (i) now
follows from another use of descending induction.

(ii) Since the lexicographical order we are using on $L[X,Y,Z]$ is
an elimination order with $X$ greater than $Y$ and $Z$, it follows
from (the analogue of) \cite[Theorem 3.6.6]{AL} that the
intersection of the set of columns of $T_d'$ with $L[Y,Z]^{{d-1}
\choose {2}}$ provides a Gr\"obner basis for $\Ima T_d \cap
L[Y,Z]^{{d-1} \choose {2}}$.

Therefore the columns of
$$
\left[ \begin{array}{c|ccccc}
  A_{d-2}  & G_{d-2} & 0 & \dots &  & 0\\
  0 & 0 & G_{d-3} & 0 & \dots & 0\\
  0 & 0 & 0 & G_{d-4} &   & 0 \\
  \vdots  & \vdots & \vdots   & & \ddots & \vdots  \\
  0  & 0 & 0 & \dots & 0 & G_1
\end{array} \right],
$$
form a Gr\"obner basis for $\Ima T_d \cap L[Y,Z]^{{d-1} \choose
{2}}$, and the claim in part (ii) follows from this.
\end{proof}

The following lemma provides motivation for part (ii) of the above
theorem.

\begin{lem}\label{se.4}
Consider the matrices $T_d$ of\/ {\rm \ref{se.1a}} and $H_d$ of\/
{\rm \ref{se.3}}.  Let $r \in L \setminus \{0\}$. Then $r$
annihilates a non-zero element of $\Coker T_d$ if and only if $r$
annihilates a non-zero element of the quotient $L[Y,Z]$-module
$\Coker H_d$ of $L[Y,Z]^{{d-1} \choose {2}}$.
\end{lem}

\begin{proof}
Suppose that $r$ annihilates a non-zero element of $\Coker T_d$.
Thus there exists a $v \in R_0^{{d-1} \choose {2}} \setminus \Ima
T_d$ such that $rv \in \Ima T_d$. We can and do assume that $v$
has been chosen so that its leading term is minimal among the
leading terms of all possible such columns. But $X  \mathbf{e}_1,
X  \mathbf{e}_2, \dots, X  \mathbf{e}_{{{d-1} \choose {2}}}$ are
all leading terms of columns of $T_d$, and so $v$ does not involve
$X$. In view of this and the fact, established in \ref{se.3}, that
$\Ima H_d \subseteq \Ima T_d$, we have $v \in L[Y,Z]^{{d-1}
\choose {2}} \setminus \Ima H_d$. Furthermore, $rv \in \Ima T_d
\cap L[Y,Z]^{{d-1} \choose {2}}$, and, by \ref{se.3}, this is the
$L[Y,Z]$-submodule of $L[Y,Z]^{{d-1} \choose {2}}$ generated by
the columns of $H_d$.

The converse is even easier.
\end{proof}

\begin{prop}\label{se.5}
For each integer $i =0, \ldots, n-1$,
$A_{i+1} A_{i+2} \dots A_{n}$ is the $(i+1)\times (n+1)$ matrix
$$\left[
\begin{array}{cccccccccc}
  Z^{n-i} & \dots & \dots & {{n-i} \choose {j}} Z^{n-i-j} Y^j & \dots & \dots
  &  Y^{n-i} &   0 & 0 & \dots \\
  0 & Z^{n-i} & \dots & & {{n-i} \choose {j}} Z^{n-i-j} Y^j &
  &  \dots & Y^{n-i} &   0 & \dots \\
   & & \ddots & & & \ddots  &  & & \ddots & \\
  \dots & 0 & 0 & Z^{n-i} & \dots & \dots &
  {{n-i} \choose {j}} Z^{n-i-j} Y^j & \dots  &
  \dots & Y^{n-i}
\end{array}\right]\! .$$
\end{prop}

\begin{proof}
The result follows from an easy reverse induction on $i$.
\end{proof}

The particular case of \ref{se.5} in which $i = 0$ yields the following.

\begin{cor}\label{se.6}
$A_{1} A_{i+2} \dots A_{n}$ is the $1\times (n+1)$ matrix whose
$(1,i+1)$-th entry is ${{n} \choose{i}} Y^i Z^{n-i}$ for all $i =
0, \ldots,n$.
\end{cor}

\begin{prop}\label{se.7}
Let $r,k \in \N$, and let $Q_{r,r+k}$ be the $r \times (r+k)$
matrix with entries in $L[Y,Z]$ given by
$$Q_{r,r+k} := \left[
\begin{array}{cccccccccc}
  Z^{k} & \dots & \dots & {{k} \choose {j}} Z^{k-j} Y^j & \dots & \dots
  &  Y^{k} &   0 & 0 & \dots \\
  0 & Z^{k} & \dots & & {{k} \choose {j}} Z^{k-j} Y^j &
  &  \dots & Y^{k} &   0 & \dots \\
   & & \ddots & & & \ddots  &  & & \ddots & \\
  \dots & 0 & 0 & Z^{k} & \dots & \dots &
  {{k} \choose {j}} Z^{k-j} Y^j & \dots  &
  \dots & Y^{k}
\end{array}\right],$$
let $\mathbf{c}_j$ denote the $j$-th column of $Q_{r,r+k}$ (for $j
= 1, \ldots, r+k$), and let $\widetilde{Q}_{r,r+k}$ be the result
of evaluation of $Q_{r,r+k}$ at $Y = Z = 1$. Thus
$$\widetilde{Q}_{r,r+k} := \left[
\begin{array}{cccccccccc}
  1 & \dots & \dots & {{k} \choose {j}}  & \dots & \dots
  &  1 &   0 & 0 & \dots \\
  0 & 1 & \dots & & {{k} \choose {j}}  &
  &  \dots & 1 &   0 & \dots \\
   & & \ddots & & & \ddots  &  & & \ddots & \\
  \dots & 0 & 0 & 1 & \dots & \dots &
  {{k} \choose {j}} & \dots  &
  \dots & 1
\end{array}\right].$$

Consider $L[Y,Z]$ as an $\nn^2$-graded ring in which
$L[Y,Z]_{(0,0)} = L$ and $\deg Y^iZ^j = (i+j,i)$. Turn the free
$L[Y,Z]$-module
$$
L[Y,Z]^r = L[Y,Z]\mathbf{e}_1 \oplus \cdots \oplus
L[Y,Z]\mathbf{e}_r
$$
into an $\nn^2$-graded module over the $\nn^2$-graded ring
$L[Y,Z]$ in such a way that $\deg \mathbf{e}_i = (0,i)$ for $i =1,
\ldots, r$. All references to gradings in the rest of this
proposition and its proof refer to this $\nn^2$-grading.

\begin{enumerate}
\item For all $i \in \nn$ and $j \in \N$, the component $\left(L[Y,Z]^r\right) _{(i,j)}$
is a free $L$-module with base
$$
\left(Y^{j-\rho}Z^{i-j+\rho}\mathbf{e}_{\rho}\right)_{\rho=\max\{j-i,1\}
,\ldots,\min \{j,r\}}.
$$
(Of course, we interpret a free module with an empty base as $0$.)
\item $\Ima Q_{r,r+k}$ is a graded submodule of $L[Y,Z]^r$, and, for
all $i \in \nn$, $j \in \N$,
$$
(\Ima Q_{r,r+k})_{(i,j)} = \begin{cases}
0 & \text{if $i < k$}, \\
\sum_{\sigma=\max\{j+k-i,1\}}^{\min \{j,r+k\}}
L Y^{j-\sigma}Z^{i-j+\sigma-k}\mathbf{c}_{\sigma} & \text{if $i \geq k$}, \\
\left(L[Y,Z]^r\right) _{(i,j)} & \text{if $i \geq 2k+r$ or $j \geq
k+r$}.
\end{cases}
$$
\item The $\nn^2$-graded $L[Y,Z]$-module $\Coker Q_{r,r+k}$
vanishes in all except finitely many degrees; in fact, $\Coker
Q_{r,r+k}$ is a finitely generated $L$-module with
$$
\Coker Q_{r,r+k} =
\bigoplus_{i=0}^{2k+r-1}\bigoplus_{j=1}^{k+r-1} \left( \Coker
Q_{r,r+k}\right) _{(i,j)}\text{;}
$$
for $0 \leq i < k$ (and $j \in \N$), the component $\left( \Coker
Q_{r,r+k}\right) _{(i,j)}$ is a free $L$-module; and for $k \leq i
\leq 2k+r-1$ and $1 \leq j \leq k+r-1$, the component $\left(
\Coker Q_{r,r+k}\right) _{(i,j)}$, as an $L$-module, is isomorphic
to the cokernel of a submatrix of $\widetilde{Q}_{r,r+k}$ made up
of the (consecutive) columns of that matrix numbered
$$\max\{j+k-i,1\},\max\{j+k-i,1\}+1, \ldots, \min \{j,r+k\}.$$
\end{enumerate}
\end{prop}

\begin{proof} (i) This is immediate from the fact that, for $\alpha, \beta
\in \nn$ and $\rho \in \{1, \ldots, r\}$, we have
$\deg Y^{\alpha}Z^{\beta}\mathbf{e}_{\rho} = (\alpha + \beta, \alpha + \rho)$.

(ii) Note that $\mathbf{c}_j$ is a homogeneous element of
$L[Y,Z]^r$ of degree $(k,j)$ (for all $j = 1, \ldots, k+r$). Hence
$\Ima Q_{r,r+k}$ is a graded submodule of $L[Y,Z]^r$, and a
homogeneous element of $\Ima Q_{r,r+k}$ is expressible as a
$L[Y,Z]$-linear combination of the columns of $Q_{r,r+k}$ in which
all the coefficients are homogeneous. Note that $\deg
Y^{\alpha}Z^{\beta}\mathbf{c}_{\sigma} = (\alpha + \beta + k,
\alpha + \sigma)$ (for $\alpha, \beta \in \nn$ and $\sigma \in
\{1, \ldots, r\}$). Hence $(\Ima Q_{r,r+k})_{(i,j)} = 0$ if $i <
k$, while $(\Ima Q_{r,r+k})_{(i,j)} =
\sum_{\sigma=\max\{j+k-i,1\}}^{\min \{j,r+k\}}L
Y^{j-\sigma}Z^{i-j+\sigma-k}\mathbf{c}_{\sigma}$ if $i \geq k$.

Notice that the vectors $Z^{k} \mathbf{e}_1, Z^{k+1} \mathbf{e}_2,
\dots, Z^{k+r-1} \mathbf{e}_{r}$ and $Y^{k+r-1} \mathbf{e}_1,
Y^{k+r-2} \mathbf{e}_2, \dots, Y^{k} \mathbf{e}_{r}$ are all in
$\Ima Q_{r,r+k}$: for any $1< s \leq r$ multiply the $s$-th
column of $Q_{r,r+k}$ by $Z^{s-1}$ and reduce with respect to
$Z^{k} \mathbf{e}_1, Z^{k+1} \mathbf{e}_2, \dots,
Z^{k+s-2}\mathbf{e}_{s-1}$ to obtain $Z^{k+s-1} \mathbf{e}_{s} \in
\Ima Q_{r,r+k}$; a similar argument shows that $Y^{k}
\mathbf{e}_{r}, Y^{k+1} \mathbf{e}_{r-1}, \dots,
Y^{k+r-1}\mathbf{e}_1 \in \Ima Q_{r,r+k}.$

Therefore, if $i \geq 2k+r$ or $j \geq k+r$, then
$$
Y^{j-\rho}Z^{i-j+\rho}\mathbf{e}_{\rho} \in \Ima Q_{r,r+k}
\quad \text{~for all~} \rho = \max\{j-i,1\} ,\ldots,\min \{j,r\}
$$
since $Y^{j-\rho}\mathbf{e}_{\rho} \in \Ima Q_{r,r+k}$ if $j \geq
k+r$, while if $j < k+r$ and $i \geq 2k+r$, then $i-j \geq k+1$
and $Z^{i-j+\rho}\mathbf{e}_{\rho} \in \Ima Q_{r,r+k}$. Thus, for
$i \geq 2k+r$ or $j \geq k+r$, all the members of the base found
in part (i) for the free $L$-module $L[Y,Z]^r _{(i,j)}$ lie in
$\Ima Q_{r,r+k}$.

(iii) All the claims of this except the final one now follow from
the previous parts. To deal with the final one, suppose that $k
\leq i \leq 2k+r-1$ and $1 \leq j \leq k+r-1$. We shall use
`overlines' to denote natural images in cokernels of elements of
free modules. It will be convenient to abbreviate $\min \{j,r+k\}$
by $\gamma$ and $\max\{j+k-i,1\}$ by $\beta$; the conditions
imposed on $i$ and $j$ ensure that $\beta \leq \gamma$.

By part (i), $\left( \Coker Q_{r,r+k}\right) _{(i,j)}$ is
generated, as $\Z$-module, by
$$
\left\{\overline{ Y^{j-\rho} Z^{i-j+\rho}\mathbf{e}_{\rho}} :
\rho= \max\{j-i,1\},\ldots, \min \{j,r\}\right\}.
$$
(Note that, once again, the conditions imposed on $i$ and $j$
ensure that $\max\{j-i,1\} \leq \min \{j,r\}$.) The fact that, for
each $\sigma=\beta, \ldots,\gamma$, we have
$\overline{Y^{j-\sigma}Z^{i-j+\sigma-k}\mathbf{c}_{\sigma}} = 0$
shows that the $\sigma$-th column of $\widetilde{Q}_{r,r+k}$ leads
to a column of relations on the generators displayed above.
Furthermore, part (ii) shows that {\em every\/} column of
relations on those generators is an $L$-linear combination of the
columns of relations arising (in this way) from the $\beta$-th,
$(\beta+1)$-th, $\ldots$, $\gamma$-th columns of
$\widetilde{Q}_{r,r+k}$.
\end{proof}

\begin{rmk}\label{se.8} Note that, with the notation of
\ref{se.7} (and provided $r > 1$), we have
$$
Q_{r-1,r+k} = Q_{r-1,r}Q_{r,r+k}.
$$
\end{rmk}

\begin{rmk}\label{se.11}
Let $B$ be a matrix with integer entries and positive rank $d$;
let $p$ be a prime number. Then it follows from the theory of the
Smith normal form that $p\Z \in \Ass_{\Z}(\Coker B)$ if and only
if the ideal generated by the $d \times d$ minors of $B$ is
contained in $p\Z$.
\end{rmk}

In view of Remark \ref{se.11} and Proposition \ref{se.7}, we are
going, in the case when $L = \Z$, to be interested in the value of
the determinant of a square matrix (with integer entries) of the
form
$$ \Omega :=
\left[\begin{array}{cccc}
  {{k} \choose{i}} & {{k} \choose{i+1}} & \dots & {{k} \choose{i+s-1}} \\
  {{k} \choose{i-1}} & {{k} \choose{i}} & \dots & {{k} \choose{i+s-2}} \\
  \vdots &  \vdots & & \vdots  \\
  {{k} \choose{i-s+1}} & {{k} \choose{i-s+2}} & \dots & {{k} \choose{i}}
\end{array}\right],
$$
where $k,s \in \N$, $i \in \nn$ and we use the convention that a
binomial coefficient ${{\xi} \choose{\eta}}$ is $0$ if either
$\eta < 0$ or $\eta > \xi$. The value of this determinant was
known to V. van Zeipel in 1865 \cite{Zeipel}; the calculation is
described in \cite[Chapter XX]{Muir}. For the convenience of the
reader, we indicate a route to the answer.

\begin{prop}\label{se.15} {\rm (See van Zeipel \cite{Zeipel}.)} Let $\Omega$ be as displayed
above. Then
$$
\det \Omega = \prod_{j=0}^{s-1} \frac{{{k+s-1-j}
\choose{i}}}{{{i+j} \choose{i}}}.
$$
\end{prop}

\begin{proof}
Add the penultimate row of $\Omega$ to the last row; in the
result, add the $(r-2)$-th row to the $(r-1)$-th, and continue in
this way until the first row has been added to the second. In this
way one sees that
$$
\det \Omega = \left|\begin{array}{cccc}
  {{k} \choose{i}} & {{k} \choose{i+1}} & \dots & {{k} \choose{i+s-1}} \\
  {{k+1} \choose{i}} & {{k+1} \choose{i+1}} & \dots & {{k+1} \choose{i+s-1}} \\
  {{k+1} \choose{i-1}} & {{k+1} \choose{i}} & \dots & {{k+1} \choose{i+s-2}} \\
  \vdots &  \vdots & & \vdots  \\
  {{k+1} \choose{i-s+2}} & {{k+1} \choose{i-s+3}} & \dots & {{k+1} \choose{i+1}}
\end{array}\right|.
$$
Now repeat the same sequence of elementary row operations, except
that, this time, stop after the second row has been added to the
third; then do a further such sequence, this time stopping after
the third row has been added to the fourth.  Continue in this way
to see that
$$
\det \Omega = \left|\begin{array}{cccc}
  {{k} \choose{i}} & {{k} \choose{i+1}} & \dots & {{k} \choose{i+s-1}} \\
  {{k+1} \choose{i}} & {{k+1} \choose{i+1}} & \dots & {{k+1} \choose{i+s-1}} \\
  {{k+2} \choose{i}} & {{k+2} \choose{i+1}} & \dots & {{k+2} \choose{i+s-1}} \\
  \vdots &  \vdots & & \vdots  \\
  {{k+s-1} \choose{i}} & {{k+s-1} \choose{i+1}} & \dots & {{k+s-1} \choose{i+s-1}}
\end{array}\right|.
$$
We proceed by induction on $i$.
When $i = 0$, it is clear from the
initial form of $\Omega$ that $\det \Omega = 1$, and the claim is true.
We therefore assume that $i > 0$, and make the obvious inductive assumption.

With reference to the last display, take out a factor $k+j-1$ from
the $j$-th row $(1 \leq j \leq s)$ and a factor $1/(i+l-1)$ from
the $l$-th column $(1 \leq l \leq s)$ to see that
$$
\det \Omega  = \frac{k(k+1)\dots(k+s-1)}{i(i+1)\dots(i+s-1)}
\left|\begin{array}{cccc}
  {{k-1} \choose{i-1}} & {{k-1} \choose{i}} & \dots & {{k-1} \choose{i+s-2}} \\
  {{k} \choose{i-1}} & {{k} \choose{i}} & \dots & {{k} \choose{i+s-2}} \\
  {{k+1} \choose{i-1}} & {{k+1} \choose{i}} & \dots & {{k+1} \choose{i+s-2}} \\
  \vdots &  \vdots & & \vdots  \\
  {{k+s-2} \choose{i-1}} & {{k+s-2} \choose{i}} & \dots & {{k+s-2} \choose{i+s-2}}
\end{array}\right| .
$$
Now use the inductive hypothesis.
\end{proof}

\begin{cor}\label{se.17p} In the situation of Proposition\/ {\rm \ref{se.15}},
we have
$$
\det \Omega = \frac{\prod_{j=0}^{s-1}{{k+s-1} \choose{i+j}} }{
            \prod_{j=0}^{s-1}{{k+s-1} \choose{j}} }.
$$
\end{cor}

\begin{proof} First note that, for $j \in \{1,\dots,s-1\}$, we have
$$
{{k+s-1-j} \choose{i}} = {{k+s-1} \choose{i+j}}
\frac{(i+1)\dots(i+j)}{(k+s-j)\dots(k+s-1)}.
$$
It therefore follows from Proposition \ref{se.15} that
\begin{align*}
\det \Omega & =  \prod_{j=0}^{s-1} \frac{{{k+s-1-j}
\choose{i}}}{{{i+j} \choose{i}}}  \\
& = \left( \prod_{j=0}^{s-1} {{k+s-1} \choose{i+j}} \frac{i!j!}{i!} \right)
\left( \prod_{j=0}^{s-1} \frac{1}{(k+s-j)\dots(k+s-1)}\right)\\
& = \left( \prod_{j=0}^{s-1} {{k+s-1} \choose{i+j}} \right)
\left( \prod_{j=0}^{s-1} \frac{j!(k+s-j-1)!}{(k+s-1)!}\right)  \\
& = \frac{\prod_{j=0}^{s-1}{{k+s-1} \choose{i+j}} }{
            \prod_{j=0}^{s-1}{{k+s-1} \choose{j}} }.
\end{align*}
\end{proof}

\begin{ntn}\label{se.16} {\rm For each $n \in \N$, we set
$$
\Pi (n) := \left\{ p : p \mbox{~is a prime factor of~}
{{n} \choose {i}} \mbox{~for some~} i \in \{0, \dots, n\}\right\}.
$$}
\end{ntn}

\begin{cor}\label{se.17} Let $r,k \in \N$ and consider the matrix
$\widetilde{Q}_{r,r+k}$ of\/ {\rm \ref{se.7}}. Let $\Delta$ be a
submatrix of $\widetilde{Q}_{r,r+k}$ formed by $c~(> 0)$
consecutive columns of that matrix; set $s := \min\{c,r\}$. If $p
\in \Z$ is a prime number such that every $s \times s$ minor of
$\Delta$ is contained in $p\Z$, then $p \in \Pi(r+k-1)$.
\end{cor}

\begin{proof} We argue by induction on $r$. Note that $\widetilde{Q}_{1,1+k}$ is
the $1 \times (1+k)$ matrix
$$ \left[
\begin{array}{cccccc}
1 & {{k} \choose {1}} & \dots & {{k} \choose {j}} & \dots & 1
\end{array} \right],
$$
and so the result is clear in this case.

Now suppose that $r >1$ and that the result has been established, for
all values of $k$, for smaller values of $r$.

If $s = r$, then there is an $r \times r$ submatrix of
$\widetilde{Q}_{r,r+k}$ of the form
$$
\Omega :=
\left[\begin{array}{cccc}
  {{k} \choose{i}} & {{k} \choose{i+1}} & \dots & {{k} \choose{i+r-1}} \\
  {{k} \choose{i-1}} & {{k} \choose{i}} & \dots & {{k} \choose{i+r-2}} \\
  \vdots &  \vdots & & \vdots  \\
  {{k} \choose{i-r+1}} & {{k} \choose{i-r+2}} & \dots & {{k} \choose{i}}
\end{array}\right],
$$
where $i \in \{0, \dots,k\}$,
such that $\det \Omega \in p\Z$. It now follows from
Corollary \ref{se.17p} that $p$ is a factor of ${{k+r-1} \choose{l}}$
for some $l \in \{0, \ldots, k+r-1\}$.

Now suppose that $s = c < r$. Set $D':= \widetilde{Q}_{r-1,r}D$.
As $\widetilde{Q}_{r-1,r+k} =
\widetilde{Q}_{r-1,r}\widetilde{Q}_{r,r+k}$ by \ref{se.8}, it
follows that $\Delta' := \widetilde{Q}_{r-1,r}\Delta$ is the
$(r-1) \times c$ submatrix of $\widetilde{Q}_{r-1,r+k}$ involving
the same columns as $\Delta$. But
$$ \widetilde{Q}_{r-1,r} = \left[\begin{array}{ccccc}
  1 & 1 & 0 & \dots & 0 \\
  0 & 1 & 1 & 0 & \dots \\
  \vdots &  & \ddots & \ddots &  \\
  0 & \dots & 0 & 1 & 1
\end{array}\right],
$$
and so the rows of $\Delta'$ are the sums of consecutive rows of
$\Delta$. Therefore any $s \times s$ minor of $\Delta'$ is the sum
of $2^s$ determinants, each one being either obviously zero or an
$s \times s$ minor of $\Delta$. Hence every $s \times s$ minor of
$\Delta'$ is contained in $p\Z$, and so, by the inductive
hypothesis, $p \in \Pi(r-1+k+1-1)$, that is, $p \in \Pi(r+k-1)$.
\end{proof}

\begin{lem}\label{se.21} The set of integers $\left\{ \#\Pi(n): n
\in \N\right\}$ is unbounded.
\end{lem}

\begin{proof} Let $(p_n)_{n \in \N}$ be an enumeration of the prime
numbers. Then, for each $n \in \N$, we have $$p_1p_2\dots p_n =
{{p_1p_2\dots p_n} \choose {1}} \in \Pi(p_1p_2\dots p_n).$$
\end{proof}

\begin{lem}\label{se.22}
Let $p \in \Z$ be a prime number. Then the sets
$$ \left\{ j \in \N : j \geq 3 \mbox{~and~} p \in \Pi(j-2)\right\}
\quad \mbox{~and~} \quad \left\{ j \in \N : j \geq 3 \mbox{~and~}
p \not\in \Pi(j-2)\right\}$$ are both infinite.
\end{lem}

\begin{proof} If $p$ divides $j-2 \in \N$, then $p\in
\Pi(j-2)$ because $p$ divides ${j-2 \choose{1}}=j-2$; hence the
first set is infinite.

To prove that the second set is infinite it is enough to show that
$p\notin \Pi(p^k-1)$ for all $k\geq 1$. Let $T$ be an
indeterminate; working modulo $p$ we have $(1+T)^{p^k-1} (1+T) =
(1+T)^{p^k} \equiv 1+T^{p^k}$ and if we compare the coefficients
of $T^{i}$ on both sides of this congruence we see that, for $0 <
i \leq p^k-1$,
$${p^k-1 \choose i} + {p^k-1 \choose i-1} \equiv 0$$
and since $p$ does not divide ${p^k-1 \choose 0}= 1 $, an easy
induction on $i$ shows that $p$ does not divide ${p^k-1 \choose
i}$ for all $i$ with $0\leq i \leq p^k-1$.
\end{proof}

We are now ready to present our main results about Singh's
example.

\begin{thm}\label{se.20} Let $R'$ denote the ring $\Z[X,Y,Z,U,V,W]/(XU + YV +
ZW)$ (considered by Singh) graded in the manner described in\/
{\rm \ref{se.1}}; let $-d \in \Z$ with $d \geq 3$; and let $p \in
\Z$ be a prime number. Then
\begin{enumerate}
\item $p\Z \in \Ass_{\Z}(H^3_{R'_+}(R')_{-d})$ if and
only if $p \in \Pi (d-2)$;
\item $\Ass_{R'_0}(H^3_{R'_+}(R')_{-d}) = \left\{(X,Y,Z)\right\}
\cup
\left\{(q,X,Y,Z) : q \in \Pi (d-2) \right\}$;
\item the set of integers
$ \left\{ \#\left(\Ass_{R'_0}(H^3_{R'_+}(R')_{-j})\right) : j \geq
3 \right\} $ is unbounded;
\item the sets
$$
\left\{ j \in \Z : j \geq 3 \mbox{~and~} (p,X,Y,Z) \in
\Ass_{R_0}(H^3_{R'_+}(R')_{-j}) \right\}
$$ and
$$
\left\{ j \in \Z : j \geq 3 \mbox{~and~} (p,X,Y,Z) \not\in
\Ass_{R_0}(H^3_{R'_+}(R')_{-j}) \right\}
$$
are both infinite; and
\item $\Ass_{R'_0}(H^3_{R'_+}(R')_{n})$ is not
asymptotically increasing for $n \rightarrow - \infty$.
\end{enumerate}
\end{thm}

\begin{proof} (i) It follows from Lemma \ref{se.1a} that
$p\Z \in \Ass_{\Z}(H^3_{R'_+}(R')_{-d})$ if and only if $p\Z \in
\Ass_{\Z}\left(\Coker T_d\right)$; furthermore, by Lemma
\ref{se.4}, this is the case if and only if $p\Z \in
\Ass_{\Z}\left(\Coker H_d\right)$, where the matrix $H_d$ is as
defined in Theorem \ref{se.3}.  It therefore follows from
\ref{se.5} and the notation introduced in \ref{se.7} that $p\Z \in
\Ass_{\Z}(H^3_{R'_+}(R')_{-d})$ if and only if
$$
p\Z \in \bigcup_{i=1}^{d-2}\Ass_{\Z} \left(\Coker Q_{i,d-1} \right).
$$

Suppose that $p \in \Pi(d-2)$, so that there exists $j \in \{1,\ldots,d-3\}$
such that $p$ is a factor of ${{d-2} \choose {j}}$. Then it follows from
Theorem \ref{se.7}(iii) that (for example) $p\Z \in \Ass_{\Z}
\left(\Coker Q_{1,d-1} \right)_{d-2,j+1}$.

Conversely, suppose that $p\Z \in \Ass_{\Z} \left(\Coker Q_{i,d-1}
\right)$, where $i \in \{1, \dots, d-2\}$. We use Theorem
\ref{se.7}(iii) to see that $p\Z \in \Ass_{\Z}(\Coker \Delta)$,
where $\Delta$ is a submatrix of $\widetilde{Q}_{i,d-1}$ formed by
$c~(> 0)$ consecutive columns of that matrix; set $s :=
\min\{c,i\}$. It follows from Proposition \ref{se.15} that
$\Delta$ has rank $s$, and therefore from Remark \ref{se.11} that
the ideal generated by the $s \times s$ minors of $\Delta$ is
contained in $p\Z$. Therefore $p \in \Pi(d-2)$ by Corollary
\ref{se.17}.

(ii) This is now immediate from part (i) and Lemma \ref{se.1a}.

(iii) This is a consequence of part (ii) and Lemma \ref{se.21}.

(iv)  This is now immediate from part (ii) and Lemmas \ref{se.1a}
and \ref{se.22}.

(v) This is a consequence of parts (ii) and (iv).
\end{proof}


\begin{thebibliography}{B-M-M}
\bibitem[A-L]{AL}W. W. Adams and P. Loustaunau, \textit{An introduction
 to Gr\"obner bases}, American Mathematical Society, Providence,
 Rhode Island, 1994.
\bibitem[B]{Brodm82}M. Brodmann, \textit{A lifting result for local
 cohomology of graded modules}, Math.\ Proc.\ Cambridge Philos.\ Soc.\
 \textbf{92} (1982), 221--229.
\bibitem[B-H]{BroHel00}M. Brodmann and M. Hellus, \textit{Cohomological patterns
 of coherent sheaves over projective schemes}, J. Pure and Appl.\ Algebra,
 to appear.
\bibitem[B-M-M]{BrMaMi00}M. Brodmann, C. Matteotti and N. D. Minh,
 \textit{Bounds for cohomological Hilbert functions of projective schemes
 over Artinian rings}, Vietnam J. Math.\ \textbf{28} (2000), 341--380.
\bibitem[B-S]{BS}M. P. Brodmann and R. Y. Sharp,
 \textit{Local cohomology:
 an algebraic introduction with geometric applications}, Cambridge
    University Press, 1998.
\bibitem[G]{Groth68}A. Grothendieck, \textit{Cohomologie locale des faisceaux
 coh\'{e}rents et th\'{e}or\`{e}mes de Lefschetz locaux et globaux (SGA 2)},
 S\'{e}minaire de G\'{e}om\'{e}trie Alg\'{e}brique du Bois-Marie 1962,
 North-Holland, Amsterdam, 1968.
\bibitem[Ma]{HMold}H. Matsumura, \textit{Commutative algebra},
 Benjamin, New York, 1970.
\bibitem[Mu]{Muir}T. Muir, \textit{The theory of determinants in the
 historical order of development}, Volume III, Macmillan, London,
 1920
\bibitem[S]{70}R. Y. Sharp, \textit{Bass numbers in the graded case, $a$-invariant
 formulas, and an analogue of Faltings' Annihilator Theorem}, J. Algebra
 \textbf{222} (1999), 246--270.
\bibitem[Si]{Singh00}A. K. Singh, \textit{$p$-torsion elements in local
 cohomology modules}, Math.\ Research Letters \textbf{7} (2000) 165--176.
\bibitem[Z]{Zeipel}V. van Zeipel, \textit{Om determinanter, hvars
 elementer \"aro binomialkoefficienter}, Lunds Universitet
 \AA rsskrift \textbf{ii} (1865) 1-68.
\end{thebibliography}
\end{document}